\documentclass[fleqn,reqno,11pt,a4paper,final]{amsart}
\usepackage[a4paper,left=35mm,right=35mm,top=30mm,bottom=30mm,marginpar=25mm]{geometry} 
\usepackage{amsmath}
\usepackage{amssymb}
\usepackage{amsthm}
\usepackage{amscd}
\usepackage[ansinew]{inputenc}
\usepackage{cite}
\usepackage{bbm}
\usepackage{color}
\usepackage[english=american]{csquotes}
\usepackage[final]{graphicx}
\usepackage{hyperref}
\usepackage{calc}
\usepackage{mathptmx}
\usepackage{bm}
\usepackage{enumerate}

\newcommand{\meantmp}[2]{#1\langle{#2}#1\rangle}
\newcommand{\mean}[1]{\meantmp{}{#1}}

\graphicspath{{../Pictures/}}

\numberwithin{equation}{section}

\newtheoremstyle{thmlemcorr}{10pt}{10pt}{\itshape}{}{\bfseries}{.}{10pt}{{\thmname{#1}\thmnumber{ #2}\thmnote{ (#3)}}}
\newtheoremstyle{thmlemcorr*}{10pt}{10pt}{\itshape}{}{\bfseries}{.}\newline{{\thmname{#1}\thmnumber{ #2}\thmnote{ (#3)}}}
\newtheoremstyle{remexample}{10pt}{10pt}{}{}{\bfseries}{.}{10pt}{{\thmname{#1}\thmnumber{ #2}\thmnote{ (#3)}}}

\theoremstyle{thmlemcorr}
\newtheorem{theorem}{Theorem}
\numberwithin{theorem}{section}
\newtheorem{lemma}[theorem]{Lemma}

\theoremstyle{thmlemcorr*}
\newtheorem{theorem*}{Theorem}
\newtheorem{lemma*}[theorem]{Lemma}
\newtheorem{corollary*}[theorem]{Corollary}
\newtheorem{proposition*}[theorem]{Proposition}
\newtheorem{problem*}[theorem]{Problem}
\newtheorem{conjecture*}[theorem]{Conjecture}
\newtheorem{definition*}[theorem]{Definition}

\theoremstyle{remexample}
\newtheorem{remark}[theorem]{Remark}

\newcommand{\Crm}{\mathrm{C}}

\newcommand{\Lrm}{\mathrm{L}}

\newcommand{\Wrm}{\mathrm{W}}

\newcommand{\Fcal}{\mathcal{F}}

\newcommand{\Lcal}{\mathcal{L}}

\newcommand{\Rcal}{\mathcal{R}}

\newcommand{\Wcal}{\mathcal{W}}

\newcommand{\Lbb}{\mathbb{L}}

\DeclareMathOperator{\diam}{diam}

\DeclareMathOperator{\diverg}{div}

\DeclareMathOperator{\supp}{supp}

\newcommand{\ee}{\mathrm{e}}

\newcommand{\setb}[2]{\bigl\{\, #1 \ \ \textup{\textbf{:}}\ \ #2 \,\bigr\}}

\newcommand{\norm}[1]{\|#1\|}

\newcommand{\abs}[1]{|#1|}

\newcommand{\absBB}[1]{\biggl|#1\biggr|}

\newcommand{\dpr}[1]{\langle #1 \rangle}

\newcommand{\dprb}[1]{\bigl\langle #1 \bigr\rangle}

\newcommand{\cl}[1]{\overline{#1}}

\newcommand{\dd}{\;\mathrm{d}}
\newcommand{\DD}{\mathrm{D}}
\newcommand{\N}{\mathbb{N}}
\newcommand{\R}{\mathbb{R}}

\newcommand{\ONE}{\mathbbm{1}}

\newcommand{\toweak}{\rightharpoonup}
\newcommand{\toweakstar}{\overset{*}\rightharpoonup}

\newcommand{\todown}{\downarrow}

\newcommand{\embed}{\hookrightarrow}
\newcommand{\cembed}{\overset{c}{\embed}}

\newcommand{\sbullet}{\begin{picture}(1,1)(-0.5,-2.5)\circle*{2}\end{picture}}
\newcommand{\frarg}{\,\sbullet\,}

\newcommand{\eps}{\epsilon}

\newcommand{\term}[1]{\textbf{#1}}

\newcounter{assumption}
\makeatletter
\newcommand{\nextas}[1]{%
  ~\refstepcounter{assumption}%
   \protected@write \@auxout{}{\string\newlabel{#1}{{(A\theassumption)}{\thepage}{(A\theassumption)}{#1}{}}}%
   \hypertarget{#1}{(A\theassumption)}%
}
\makeatother

\def\Xint#1{\mathchoice 
{\XXint\displaystyle\textstyle{#1}}%
{\XXint\textstyle\scriptstyle{#1}}%
{\XXint\scriptstyle\scriptscriptstyle{#1}}%
{\XXint\scriptscriptstyle\scriptscriptstyle{#1}}%
\!\int} 
\def\XXint#1#2#3{{\setbox0=\hbox{$#1{#2#3}{\int}$} 
\vcenter{\hbox{$#2#3$}}\kern-.5\wd0}} 
\def\dashint{\,\Xint-}

\renewcommand{\eps}{\varepsilon}
\renewcommand{\epsilon}{\varepsilon}
\renewcommand{\phi}{\varphi}
\renewcommand{\hat}{\widehat}

\begin{document}


\title[Existence and regularity for rate-independent systems]{Existence and regularity for a class of\\rate-independent systems}

\author{Filip Rindler}
\address{\textit{Filip Rindler:} Mathematics Institute, University of Warwick, Coventry CV4 7AL, UK.}
\email{F.Rindler@warwick.ac.uk}

\author{Sebastian Schwarzacher}
\address{\textit{Sebastian Schwarzacher:} Katedra matematick\'{e} anal\'{y}zy, Charles University Prague, Sokolovsk\'{a} 83, 186 75 Praha 8, Czech Republic.}
\email{schwarz@karlin.mff.cuni.cz}


\hypersetup{
  pdfauthor = {Filip Rindler (University of Warwick) and Sebastian Schwarzacher (Charles University)},
  pdftitle = {Existence and regularity for a class of rate-independent systems},
  pdfsubject = {},
  pdfkeywords = {}
}


\maketitle
\thispagestyle{empty}


\begin{abstract}
Despite the many applications of rate-independent systems, their regularity theory is still largely unexplored. Usually, only weak solution with potentially very low regularity are considered, which requires non-smooth techniques. In this work, however, we directly prove the existence of H\"{o}lder-regular strong solutions for a class of rate-independent systems. We also establish further assertions about higher regularity of our solutions. The proof proceeds via a time-discrete Rothe approximation, careful elliptic regularity estimates in the discrete situation and evolutionary techniques.
\vspace{4pt}



\vspace{4pt}

\noindent\textsc{Date:} \today{}. 
\end{abstract}


\section{Introduction}

Rate-independent systems are used to model a plethora of physical phenomena for which the speed of the evolution does not influence the amount of energy dissipation, including elasto-plasticity, damage \& delamination in solids, crack propagation, and shape-memory alloys. We refer to the recent monograph~\cite{MielkeRoubicek15book} for an up-to-date overview over the ample literature on both theory and applications of rate-independent systems. In this work we consider only \emph{purely dissipative} systems, i.e.\ those without elastic variables (which by definition are those quantities that can be changed without dissipating energy).

A purely formal prototype for a rate-independent system of this type is
\[
 \frac{\dot{u}}{\abs{\dot{u}}} - \Delta u + \DD W_0(u) = f, \qquad u \colon [0,T] \times \Omega \to \R^m.
\]
For the mathematical analysis, the crucial feature of such systems of equations is the quasi-static nature of their evolution, namely that solutions simply \emph{rescale} when scaling the external forces. In a sense, quasi-static evolutions only \emph{follow} the total (energetic and external) forces and hence the evolution should in some cases \enquote{inherit} regularity from the external force. Questions about the regularity of solutions were fact already raised at the very beginning of the modern energetic theory of rate-independent systems, see in particular~\cite{MielkeTheilLevitas02} and~\cite{MielkeTheil04} (Section~7.3 in the second reference discusses temporal regularity for the uniformly convex case). Further investigations in that direction were carried out recently in~\cite{MielkeZelik14}, but not much else appears to be known in general.

We here advance the existence and regularity theory for rate-independent systems by showing existence of strong solutions with essential optimal regularity (in space and time) for a class of rate-independent systems with a quadratic gradient regularizer. Our assumptions apply for instance to some \enquote{mild} double-well energies, where the range of allowed non-convexity depends on the $\Lrm^2$-Poincar\'{e} embedding constant of the domain. Our condition entails that the \emph{regularized} functional is convex. This framework is used frequently, as some regularization is already needed for the existence of solutions, see~\cite{MielkeRoubicek15book}.

While the theory of rate-independent systems is dominated by several notions of weak solution concepts (see~\cite{Mielke11} for an overview), very little is known about the existence of strong solutions. In particular, in the non-convex case it is perhaps surprising that they indeed exist. Even more surprising is that quite a lot of regularity can be established, namely H\"{o}lder continuity in space and time. Our results should open up new applications both in the theory and for numerical approximations.

For technical reasons we only consider the cases of two and three spatial dimensions. Some results are also true in higher dimensions, but additional (more restrictive) assumptions are necessary. In view of the fact that the two- and three-dimensional situations are the most physically relevant anyway, we confine ourselves to these cases.

Concretely, we will investigate the class of rate-independent system that can be formulated as follows: For $\Omega \subset \R^d$ a bounded Lipschitz domain, $d=2,3$, and $T > 0$ consider the (formal) system
\begin{align} \label{eq:PDE}
  \left\{
  \begin{aligned}
    \partial \Rcal_1(\dot{u}(t)) - \Lcal_t u(t) + \DD\Wcal_0(u(t)) &\ni f
      \quad \text{in $[0,T] \times \Omega$,}\\
    u(t)|_{\partial \Omega} &= 0  \quad\text{for $t \in [0,T]$,} \\
    \qquad u(0) &= u_0.
  \end{aligned} \right.
\end{align}
Here, $\Rcal_1 \colon \Lrm^1(\Omega;\R^m) \to \R \cup \{+\infty\}$ is the \term{rate-independent dissipation potential}, which is assumed to be proper (not identically $+\infty$), convex, and positively $1$-homogeneous; $\partial \Rcal_1$ is its subdifferential. By $\Wcal \colon \Lrm^q(\Omega;\R^m) \to \R$, $q \in (1,\infty)$, we denote the \term{(elastic) energy functional}, $f \in \Wrm^{1,a}(0,T;\Lrm^p(\Omega;\R^m))$, for $p \in [2,\infty)$ and $a \in (1,\infty)$, is the \term{external loading (force)}, and $u_0 \in (\Wrm^{1,2}_0 \cap \Lrm^q)(\Omega;\R^m)$ is the \term{initial value}. Finally, the \term{regularizer} $\Lcal_t$ is a (possibly time-dependent) second-order linear PDE operator in the space variables (most commonly, $\Lcal_t = \Delta$). Precise assumptions are detailed below. Note that the assumption of zero Dirichlet boundary values is imposed merely to simplify the exposition, analog results for non-zero Dirichlet boundary values cause only technical changes of the arguments below.

We call a map $u \in \Lrm^\infty(0,T;(\Wrm^{2,2}_0 \cap \Lrm^q)(\Omega;\R^m))$ such that its weak time derivative $\dot{u}$ has regularity $\dot{u} \in \Lrm^1(0,T;\Wrm^{1,2}(\Omega;\R^m))$ a \term{strong solution} to~\eqref{eq:PDE} if
\[
  \Lcal_t u(t) - \DD\Wcal_0(u(t)) + f(t) \in \Lrm^\infty(0,T;\Lrm^2(\Omega;\R^m))
\]
and
\[
  \left\{ \begin{aligned}
    \Lcal_t u(t) - \DD\Wcal_0(u(t)) + f(t) &\in \partial \Rcal_1(\dot{u}(t))  \qquad\text{for a.e.\ $t \in [0,T]$,} \\
    u(0) &= u_0.
  \end{aligned} \right.
\]
Writing out the definition of the subdifferential, the above inclusion means
\begin{equation} \label{eq:strongsol_ineq}
  \Rcal_1(\dot{u}(t)) + \dprb{\Lcal_t u(t) - \DD\Wcal_0(u(t)) + f(t), \xi(t) - \dot{u}(t)} \leq \Rcal_1(\xi(t))
\end{equation}
for all $\xi \in \Lrm^1(0,T;\Wrm^{1,2}_0(\Omega;\R^m))$ and almost every $t \in [0,T]$, where $\dpr{\frarg,\frarg}$ denotes the $\Lrm^2$-duality pairing. We recall that the second condition makes sense since our regularity assumptions on $u, \dot{u}$ imply that $u \in \Crm([0,T];\Lrm^2(\Omega;\R^m))$.

In the remainder of the introduction, we will state our assumptions and main results.



\subsection{Assumptions} \label{sc:assume}

The following conditions are assumed in all of the following, if not stated otherwise:

\begin{enumerate}[({A}1)]
  \item Let $d\in\{2,3\}$ and let $\Omega \subset \R^d$ be open, bounded and with boundary of class $\Crm^{1,1}$.
  \item The rate-independent dissipation (pseudo)potential $\Rcal_1 \colon \Lrm^1(\Omega;\R^m) \to \R \cup \{+\infty\}$ is given as
\[\qquad
  \Rcal_1(v) = \int_\Omega R_1(v(x)) \dd x,  \qquad v \in \Lrm^1(\Omega;\R^m),
\]
with $R_1 \colon \R^m \to \R \cup \{+\infty\}$ proper, convex, lower semicontinuous, and positively $1$-homogeneous, i.e.\ $R_1(\alpha w) = \alpha R_1(w)$ for any $\alpha \geq 0$ and $w \in \R^m$.


  \item The energy functional $\Wcal_0 \colon \Lrm^q(\Omega;\R^m) \to [0,\infty]$, where $q \in (1,\infty)$, has the form
\[\qquad
  \Wcal_0(u) = \int_\Omega W_0(u(x)) \dd x
\]
with $W_0 \in \Crm^1(\R^m;[0,\infty))$ satisfying the following assumptions for constants $C,\mu > 0$:
\begin{align}
  \qquad
  \Crm^{-1}(\abs{v}^q-1) &\leq W_0(v) \leq C(\abs{v}^q+1),   \label{eq:Wgrowth}\\
  \abs{\DD W_0(v)}  &\leq C(1+\abs{v}^{q-1}),               \label{eq:DWgrowth}\\
 -\mu \abs{v-z}^2 &\leq (\DD W_0(v)-\DD W_0(z))\cdot(v-z),  \label{eq:Wmon}\\
  \mu C_P(\Omega)^2 &< 1.                                    \label{eq:muCP}
\end{align}
Here, $C_P(\Omega) > 0$ is the (best) $\Lrm^2$-Poincar\'{e} embedding constant of $\Omega$, i.e.\ the smallest $C > 0$ such that
\[\qquad
  \norm{v}_{\Lrm^2} \leq C \norm{\nabla v}_{\Lrm^2}  \qquad
  \text{for all $v \in \Wrm^{1,2}_0(\Omega;\R^m)$.}
\]
These assumptions are for instance satisfied even for the non-convex double-well potential $W_0(s) = \gamma (s^2-1)^2$ for $0 < \gamma < (2C_P(\Omega))^{-2}$ (and $q = 4$), as a straightforward calculation shows.

 \item The regularizer $\Lcal_t$ is a second-order linear PDE operator of the form
\[ \qquad
  [\Lcal_t v]^{\beta} := \sum_j \partial_j \sum_{\alpha,i} A^{\alpha, \beta}_{i,j}(t,x) \, \partial_i v^\alpha,  \qquad
  1 \leq \alpha, \beta \leq m, \;\; 1 \leq i,j \leq d.
\]
We assume the coefficients $A^{\alpha,\beta}_{i,j}$ to satisfy the following continuity, ellipticity and symmetry conditions:
\begin{align} \qquad
  A^{\alpha,\beta}_{i,j} = A^{\beta,\alpha}_{j,i} \in \Crm^{0,1}([0,T]\times\overline{\Omega})
    &\quad \text{for $\alpha,\beta\in \{1,...,N\}$, $i,j\in\{1,...,d\}$,} \label{eq:Lcoeff} \\
  \sum_{i,j,\alpha,\beta} \xi^\alpha_i A^{\alpha, \beta}_{i,j}(t,x) \xi^\beta_j \geq \kappa \abs{\xi}^2
    &\quad\text{for $\xi \in \R^{m \times d}$ and some $\kappa > 0$.}  \label{eq:Lell} 
\end{align}
  
  \item The external force has regularity $f \in \Wrm^{1,a}(0,T;\Lrm^p(\Omega;\R^m))$ for $a \in (1,\infty)$, $p \in [2,\infty)$.
  
  \item The initial value satisfies $u_0 \in (\Wrm^{1,2}_0 \cap \Lrm^q)(\Omega;\R^m)$.
\end{enumerate}

\begin{remark}
Note that we are allowing non-convexity in $W_0$, but by~\eqref{eq:Wmon}~\eqref{eq:muCP}, this non-convexity cannot be too strong. It can be shown that our conditions entail that the combined energy functional $\Wcal(u) := \frac{1}{2} \norm{\nabla u}_{\Lrm^2} + \Wcal_0(u)$ is convex. Some convexity is also necessary, since for strongly non-convex $W_0$ counterexamples to regularity exist, see~\cite{MielkeTheil04}.
\end{remark}

\subsection{Main result}

In the course of this work we will prove the following result about the existence and regularity of a strong solution to~\eqref{eq:PDE}:

\begin{theorem} \label{thm:main}
Under the above assumptions, there exists a strong solution
\[
  u \in \Lrm^\infty(0,T;(\Wrm^{2,2}_0 \cap \Lrm^q)(\Omega;\R^m))
  \qquad\text{with}\qquad
  \dot{u} \in \Lrm^1(0,T;\Wrm^{1,2}(\Omega;\R^m))
\]
to~\eqref{eq:PDE}. Moreover, this solution has the following additional regularity properties:
\begin{enumerate}[(i)]
\item $\nabla^2 u\in \Lrm^\infty (0,T;\Lrm^p(\Omega;\R^m))$,
\item $\nabla \dot{u}\in \Lrm^a(0,T;\Lrm^2(\Omega;\R^m))$,
\item $u\in \Crm^{0,\gamma}([0,T]\times \overline{\Omega};\R^m)$ for some $\gamma \in (0,1)$,
\item If $p>d$, then $\nabla u\in \Crm^{0,\zeta}([0,T]\times \overline{\Omega};\R^m)$ for some $\zeta \in (0,1)$.
\end{enumerate}
\end{theorem}


\begin{remark}
Our proof also provides the following quantitative estimates:
\begin{align*}
\norm{\nabla^2 u}_{\Lrm^\infty(\Lrm^p)}&\leq C\Bigl(1+\norm{f}_{\Lrm^\infty(\Lrm^p)}+\norm{f}_{\Lrm^\infty(\Lrm^2)}^{q-1}\Bigr),\\
\norm{\nabla \dot{u}}_{\Lrm^a(\Lrm^2)}&\leq C\Bigl(1+\norm{f}_{\Wrm^{1,a}(\Lrm^2)}\Bigr).
\end{align*}
Here, the constant $C>0$ depends on all constants in the assumptions and on $p, a, T, \abs{\Omega}$. The estimates follow from the discrete estimates~\eqref{eq:apri_space},~\eqref{eq:apri_time} below.
 
The oscillation estimates are quantified in the following manner:
\[
[u]_{\Crm^{0,\gamma}([0,T]\times\overline{\Omega})}\leq C\Big(1 + \norm{u}_{\Lrm^\infty(\Wrm^{2,p})}+\norm{\dot{u}}_{\Lrm^a(\Wrm^{1,2})}\Big),
\]
where $\gamma \in (0,1)$, and, if $p>d$,
\[
[\nabla u]_{\Crm^{0,\zeta}([0,T]\times\overline{\Omega})}\leq C\Big(\norm{\nabla u}_{\Lrm^\infty(\Wrm^{1,p})}+\norm{\nabla\dot{u}}_{\Lrm^a(\Lrm^{2})}\Big),
\]
where $\zeta \in (0,1)$. These estimates are explained in Remark~\ref{rem:H} below.
\end{remark}

Our proof proceeds via a Rothe time-discretization scheme and crucial \enquote{elliptic} estimates at the discrete level, see Lemma~\ref{lem:apriori-space}. We give a brief formal overview over the estimates that can be expected in Section~\ref{sc:formal}, then proceed to the rigorous proof in Sections~\ref{sc:estimate} and~\ref{sc:proof}.

\subsection*{Acknowledgments} 

The authors would like to thank Florian Theil for interesting discussions related to this work. F.~R.\ gratefully acknowledges the support from an EPSRC Research Fellowship on ``Singularities in Nonlinear PDEs'' (EP/L018934/1). S.~S.\ thanks the program PRVOUK~P47 of the Charles University Prague.

\section{Formal a-priori estimates}  \label{sc:formal}

We first illustrate what can be gained from a-priori estimates by purely non-rigorous, formal arguments. For the purpose of illustration, we also restrict ourselves to the case $\Lcal_t = \Delta$. All of these arguments will be made precise in the following sections. So, assume that we have a smooth $u \colon [0,T] \times \Omega \to \R^m$ satisfying
  \begin{align}
&\int_\Omega R_1(\dot{u}(t))-\nabla u(t)\cdot \nabla(\xi-\dot{u}(t))+[-\DD W_0(u(t))+f(t)]\cdot(\xi-\dot{u}(t)) \dd x  \notag\\
&\qquad\leq \int_\Omega R_1(\xi) \dd x  \label{eq:toy}
\end{align}
for all smooth $\xi \colon \Omega \to \R^m$ and all $t \in [0,T]$.

\subsection{Estimates in space}
Concerning the space, the weak formulation~\eqref{eq:toy} implies that for almost every $(t,x) \in [0,T] \times \Omega$ there is  $z(t,x)\in \partial R_1(\dot{u}(t,x))$ such that
\[
  z(t) - \Delta u(t)+\DD W_0(u(t))=f(t).
\]
Since by the properties of $R_1$ we have $\abs{z(t,x)}\leq C$ uniformly, we get that
$- \Delta u(t)\in \Lrm^\infty$ if $f$ and $\DD W_0$ are globally bounded.
This implies that $\nabla^2 u(t)\in \Lrm^s$ for all $s\in(1,\infty)$. Moreover, for $\Omega$ bounded with smooth boundary, we have the estimate 
\[
\norm{\nabla^2u(t)}_{\Lrm^s}\leq C\norm{\Delta u(t)}_{\Lrm^s} \qquad
 \text{for $s\in (1,\infty)$,}
\]
and the constant depends only on $\Omega,s,d$, see~\cite[(11.8)]{LadyzhenskayaUraltseva68}, or~\cite{Browder60,Browder61}. 
Thus, if for the moment $\DD W_0$ is assumed bounded, we find
\begin{align*}
 \norm{\nabla^2u(t)}_{\Lrm^s}\leq C(1+\norm{f(t)}_{\Lrm^s}).
\end{align*}
The usual embedding results then yield
\[
  \nabla u \in \Lrm^\infty(0,T;\Crm^{0,\alpha}(\cl{\Omega};\R^m)) \quad \text{for all $\alpha\in[0,1)$.}
\]
The above derivation and the extension to $\DD W_0$ unbounded but satisfying~\eqref{eq:DWgrowth},~\eqref{eq:Wmon}, is made precise in Lemma~\ref{lem:apriori-space}.


\subsection{Testing with $u$.}
\label{ssec:u}
Choosing $\xi=\dot{u}(t)-\phi$ for any $\phi\in \Crm^\infty_0(\Omega)$ implies
  \begin{align*}
\int_\Omega R_1(\dot{u}(t))-R_1(\dot{u}(t)-\phi)+\nabla u(t)\cdot \nabla\phi+[\DD W_0(u(t))-f(t)]\cdot\phi \dd x\leq 0,
\end{align*}
which gives, by adding $R_1(-\phi)$ on both sides and invoking the subaditivity of $R_1$ (since it is convex and positively $1$-homogeneous), that
\begin{equation} \label{eq:toy_phi}
\int_\Omega \nabla u(t)\cdot \nabla\phi+[\DD W_0(u(t))-f(t)]\cdot\phi \dd x\leq \int_\Omega R_1(-\phi)\dd x.
\end{equation}
Using $\phi=u(t)$, we get by~\eqref{eq:DWgrowth} that
\begin{align*}
\norm{\nabla u(t)}_{\Lrm^2}^2\leq \norm{f(t)}_{\Lrm^2}\norm{u(t)}_{\Lrm^2}+C\norm{u(t)}_{\Lrm^1}
+C\norm{u(t)}_{\Lrm^2}^2,
\end{align*}
which implies by Poincar\'{e}'s inequality that
\[
\norm{\nabla u(t)}_{\Lrm^2}\leq \norm{f(t)}_{\Lrm^2}+C(1+\norm{u(t)}_{\Lrm^2}).
\]
This holds at almost every $t \in [0,T]$. Observe that in the case $W_0$ is convex the term $\DD W_0(u(t))\cdot u$ is even positive. In this case or in the case where $\DD W_0$ is a-priori bounded, we get an estimate for $\norm{\nabla u(t)}_{\Lrm^2}$ that is is independent of $u(t)$. 

\subsection{Testing with $\dot{u}$.}
A higher-order a-priori estimate can be derived by choosing $\xi=0$. Then, 
 \begin{align*}
\int_\Omega R_1(\dot{u}(t))+\partial_t\biggl(\frac{\abs{\nabla u(t)}^2}{2}\biggr)+\partial_t W_0(u(t))\dd x\leq \int_\Omega f(t)\cdot \dot{u}(t)\dd x.
\end{align*}
We integrate over the time intervall $(0,\tau) \subset (0,T)$ and, assuming for the simplicity that $u(0) = 0$ and $W_0(0)=0$, we find that
 \begin{align*}
&\int_0^\tau\!\!\int_\Omega R_1(\dot{u})\dd x\dd t+ \int_\Omega\frac{\abs{\nabla u(\tau)}^2}{2}+W_0(\tau)\dd x\\
&\qquad=\int_0^\tau\!\!\int_\Omega f\cdot
 \dot{u}\dd x\dd t\\
&\qquad= -\int_0^\tau\!\!\int_\Omega \dot{f}\cdot
 u\dd x\dd t+\int_\Omega f(\tau)\cdot u(\tau) \dd x
\\
&\qquad \leq \norm{\dot{f}}_{\Lrm^1(\Lrm^2)}\norm{u}_{\Lrm^\infty (\Lrm^2)}+\norm{f}_{\Lrm^\infty(\Lrm^2)}\norm{u}_{\Lrm^\infty (\Lrm^2)}.
\end{align*}
This implies, by taking the supremum over all $\tau \in [0,T]$, and absorbing via the Poincar\'{e} inequality, that
\begin{align*}
\norm{\dot{u}}_{\Lrm^1(\Lrm^1)}+\norm{\nabla u}^2_{\Lrm^\infty(\Lrm^2)}+\sup_{\tau}\int_\Omega W_0(u(\tau))&\leq C \Bigl( \norm{\dot{f}}_{\Lrm^1(\Lrm^2)}^2+ \norm{f}_{\Lrm^\infty(\Lrm^2)}^2 \Bigr) \\
&\leq C\norm{f}_{\Wrm^{1,1}(\Lrm^2)}^2.
\end{align*}

\subsection{Testing with $\ddot{u}$}
The following estimate is the crucial one, yet its derivation in continuous time is surprisingly simple. More effort will be needed later to derive it for the time-discretized situation. We differentiate~\eqref{eq:toy_phi} in time to find
 \begin{align*}
\int_\Omega \nabla \dot{u}(t)\cdot \nabla\phi+[\DD^2 W_0(u(t))\dot{u}(t)-\dot{f}(t)]\cdot\phi \dd x\leq 0.
\end{align*}
Now we use for a fixed $t \in [0,T]$ the test function $\phi = \dot{u}(t)$, whereby
\begin{align*}
\norm{\nabla \dot{u}(t)}_{\Lrm^2}^2+\int_\Omega \DD^2W_0(u(t))[\dot{u}(t),\dot{u}(t)] \dd x \leq \norm{\dot{f}(t)}_{\Lrm^2}\norm{\dot{u}(t)}_{\Lrm^2}.
\end{align*}
We remark that this is exactly the point where the mild convexity assumption~\eqref{eq:muCP} is essential. Indeed,~\eqref{eq:Wmon} implies (via difference quotients) in the case where $W_0$ is twice differentiable that $\DD^2W_0(u)[\dot{u},\dot{u}]\geq -\mu \abs{\dot{u}}^2$. Hence we find by H\"older's inequality and the Poincar\'{e}-Friedrich inequality that
\begin{align*}
\norm{\nabla \dot{u}(t)}_{\Lrm^2}^2 &\leq \norm{\dot{f}(t)}_{\Lrm^2}\norm{\dot{u}(t)}_{\Lrm^2} +\mu \norm{\dot{u}(t)}_{\Lrm^2}^2\\
&\leq C_P(\Omega)\norm{\dot{f}(t)}_{\Lrm^2}\norm{\dot{u}(t)}_{\Lrm^2} +\mu C_P(\Omega)^2\norm{\nabla\dot{u}(t)}_{\Lrm^2}^2
\end{align*}
Thus,
\begin{align*}
\norm{\nabla \dot{u}(t)}_{\Lrm^2}\leq \frac{C}{1-\mu C_P(\Omega)^2}\norm{\dot{f}(t)}_{\Lrm^2},
\end{align*}
and the constant is positive by~\eqref{eq:muCP}.
%
%

\section{Existence and regularity of solutions}  \label{sc:estimate}

In this section we will prove Theorem~\ref{thm:main}. We will do so by the Rothe method and discrete analogues of the estimates from the previous section. It will turn out that our a-priori information on solutions is quite strong and thus we obtain compactness in a variety of spaces. The difficulty is to establish the limit equation and for this we will need H\"{o}lder continuity of the solution.

\subsection{Time discretization}

We consider a sequence of partitions
\[
  0 = t^N_0 < t^N_1 < \cdots < t^N_N = T,  \qquad\text{where}\qquad
  t^N_k - t^N_{k-1} = \frac{T}{N}, \quad N \in \N,
\]
and look for corresponding discrete-time approximations
\[
  (u^N_k)_{k=0,\ldots,N} \subset (\Wrm^{1,2}_0 \cap \Lrm^q)(\Omega;\R^m),
\]
which solve a suitable discrete version of~\eqref{eq:PDE}. As the approximations for the external force $f$ we set
\[
  f^N_k := f(t^N_k) \quad\text{for $k = 0,\ldots,N$,}
  \qquad\text{and}\qquad
  f^N := \sum_{k=1}^N \ONE_{(t^N_{k-1},t^N_k]}f^N_k.
\]
Observe that our assumption $f \in \Wrm^{1,a}(0,T;\Lrm^p(\Omega;\R^m))$ for some $a \in (1,\infty)$ implies $f^N\in \Lrm^\infty(0,T;\Lrm^p(\Omega;\R^m))$. Further define the following approximations of the elliptic operator $\Lcal_t$:
\[
  [\Lcal^N_kv]^\beta = \diverg([\Lbb^N_kv]^\beta)
  := \sum_j \partial_j \sum_{i,\alpha} A^{\alpha,\beta}_{i,j}(t_k^N,\frarg)\partial_i v^\alpha.
\]


Now, iteratively at each $k = 0,1,\ldots,N$, minimize the functional
\[
  \Fcal^N_k(v) := \int_\Omega R(v-u^N_{k-1}) + \nabla v : \frac{\Lbb^N_k}{2} : \nabla v + W_0(v) - f^N_k\cdot v \dd x
\]
over all $v \in (\Wrm^{1,2}_0 \cap \Lrm^q)(\Omega;\R^m)$. Here, we used the notation
\[
  \nabla v : \Lbb^N_k : \nabla w := \sum_{\alpha,\beta,i,j} A^{\alpha, \beta}_{i,j}(t_k^N,x) \partial_i v^\alpha \partial_j w^\beta.
\]
Since $R$ is convex and lower semicontinuous and $W_0$ is of lower order, we may deduce by the usual Direct Method that a minimizer exists, which we call $u^N_k$. More precisely, we take a minimizing sequence $(v_j) \subset  (\Wrm^{1,2}_0 \cap \Lrm^q)(\Omega;\R^m)$ with $\Fcal^N_k(v_j) \to \min \Fcal^N_k$. Then, by the coercivity of $W_0$ (see~\eqref{eq:Wgrowth}), the strong ellipticity of $\Lcal_t$ (see~\eqref{eq:Lell}) and $R \geq 0$, we get the estimate
\[
  \norm{\nabla v_j}_{\Lrm^2}^2 + \norm{v_j}_{\Lrm^q}^q \leq C(1 + \norm{f^N_k}_{\Lrm^2} \cdot \norm{v_j}_{\Lrm^2})
\]
for a $j$-independent constanct $C > 0$. Thus, using the Poincar\'{e} and Young inequalities,
\begin{equation} \label{eq:coerc}
  \norm{\nabla v_j}_{\Lrm^2} \leq C(1+\norm{f^N_k}_{\Lrm^2})  \qquad\text{and}\qquad
  \norm{v_j}_{\Lrm^q}^q \leq C(1+\norm{f^N_k}_{\Lrm^2}^2).
\end{equation}
That is, we have shown coercivity in $(\Wrm^{1,2}_0 \cap \Lrm^q)(\Omega;\R^m)$. Hence, we may assume after selecting a non-relabelled subsequence that $v_j \toweak v$ in $\Wrm^{1,2} \cap \Lrm^q$. Furthermore, by the compact embedding $\Wrm^{1,2}_0(\Omega;\R^m) \cembed \Lrm^2(\Omega;\R^m)$ and selecting another subsequence, $v_j \to v$ pointwise almost everywhere. Now, for the convex terms in $\Fcal^N_k$, we get lower semicontinuity immediately and for $W_0$ ($\geq 0$) we estimate
\[
  \liminf_{j\to\infty} \int_\Omega W_0(v_j(x)) \dd x
  \geq \int_\Omega W_0(v(x)) \dd x,
\]
by Fatou's lemma.
Hence, the Direct Method applies and yields the existence of a minimizer, which we call $u^N_k$.

The minimizer $u^N_k$ satisfies the Euler--Lagrange equation
\[
  0 \in \partial R(u^N_k-u^N_{k-1}) - \Lcal^N_k u^N_k + \DD W_0(u^N_k) - f^N_k,
\]
in a weak sense. That is, for any test function $\xi \in \Wrm^{1,2}_0(\Omega;\R^m)$ it holds that
\begin{align}
  &\int_\Omega R(u^N_k-u^N_{k-1}) - \nabla u^N_k : \Lbb^N_k : \nabla(\xi-(u^N_k-u^N_{k-1})) \notag\\
  &\qquad + \bigl[ -\DD W_0(u^N_k) + f^N_k \bigr] \cdot (\xi-(u^N_k-u^N_{k-1})) \dd x \leq \int_\Omega R(\xi) \dd x.  \label{eq:ELweak}
\end{align}
To see this, we observe that first for $\xi \in (\Wrm^{1,2}_0 \cap \Lrm^q)(\Omega;\R^m)$ we have
\begin{equation} \label{eq:Fineq}
  0 \leq \frac{\Fcal^N_k \bigl( u^N_k + \eps(\xi + u^N_{k-1}-u^N_k) \bigr) - \Fcal^N_k(u^N_k)}{\eps},  \qquad \eps > 0.
\end{equation}
First, since $R$ is homogeneous of degree one and convex, it is subadditive, i.e.\ $R(a+b) \leq R(a)+R(b)$, and so
\begin{align*}
  &R \bigl( u^N_k + \eps(\xi + u^N_{k-1}-u^N_k) - u^N_{k-1} \bigr) - R \bigl(u^N_k - u^N_{k-1} \bigr) \\
  &\qquad = R \bigl( \eps\xi + (1-\eps)(u^N_k - u^N_{k-1}) \bigr) - R \bigl(u^N_k - u^N_{k-1} \bigr) \\
  &\qquad \leq \eps R(\xi) - \eps R(u^N_k - u^N_{k-1}).
\end{align*}
For the regularizer we may compute using the symmetry of the coefficients in $\Lcal_t$, see~\eqref{eq:Lcoeff}, and setting $\eta := \xi + u^N_{k-1}-u^N_k$,
\begin{align*}
  &\frac{1}{\eps} \int_\Omega [\nabla u^N_k + \eps \nabla \eta] : \frac{\Lbb^N_k}{2} : [\nabla u^N_k + \eps \nabla \eta] - \nabla u^N_k : \frac{\Lbb^N_k}{2} : \nabla u^N_k \dd x  \\
  &\qquad\to \int_\Omega \nabla u^N_k : \Lbb^N_k : \nabla \eta \dd x  \qquad
  \text{as $\eps \todown 0$}
\end{align*}
by the $\Lrm^2$-bounds on all involved quantities. Finally, note (we let $\dashint_0^\eps:=\frac{1}{\eps}\int_0^\eps$)
\begin{align*}
  \frac{1}{\eps} \int_\Omega W_0(u^N_k + \eps \eta) - W_0(u^N_k) \dd x
  &= \int_\Omega \dashint_0^\eps \DD W_0(u^N_k + \tau \eta) \cdot \eta \dd \tau \dd x \\
  &\to \int_\Omega \DD W_0(u^N_k) \cdot \eta \dd \tau \dd x
\end{align*}
by the continuity of $\DD W_0$ and the estimate~\eqref{eq:Wgrowth} on the growth of $\DD W_0$. Thus, letting $\eps \todown 0$ in~\eqref{eq:Fineq}, we arrive at~\eqref{eq:ELweak} for $\xi \in (\Wrm^{1,2}_0 \cap \Lrm^q)(\Omega;\R^m)$. A density argument allows us to conclude~\eqref{eq:ELweak} also for $\xi \in \Wrm^{1,2}_0(\Omega;\R^m)$.

\subsection{Discrete a-priori estimates in space} \label{ssec:disc-aprior}
Set
\[
  \delta^N_k:=\frac{u_{k}^N-u_{k-1}^N}{h},  \qquad h := t^N_k - t^N_{k-1}.
\]
Then, dividing~\eqref{eq:ELweak} by $h$ and replacing $\xi/h$ by $\xi$, we get
\begin{align}
&\int_\Omega R(\delta_k^N) - \nabla u^N_k : \Lbb^N_k : \nabla(\xi-\delta^N_k)+\bigl[ -\DD W_0(u^N_k) + f^N_k \bigr]\cdot(\xi-\delta^N_k) \dd x \notag\\
&\qquad \leq \int_\Omega R(\xi) \dd x.  \label{eq:discrete1}
\end{align}
Now, we may further replace $\xi - \delta^N_k$ by $\hat{\xi}$ and use the subadditivity of $R$ to get
\[
\int_\Omega - \nabla u^N_k : \Lbb^N_k : \nabla \hat{\xi}+\bigl[ -\DD W_0(u^N_k) + f^N_k \bigr] \cdot \hat{\xi} \dd x \leq \int_\Omega R(\hat{\xi}) \dd x.  
\]
We localize this with a cut-off function in $\hat{\xi}$ to infer
\[
- \nabla u^N_k : \Lbb^N_k : \nabla \hat{\xi}+\bigl[ -\DD W_0(u^N_k) + f^N_k \bigr] \cdot \hat{\xi} \leq R(\hat{\xi}) \leq C\abs{\hat{\xi}}  \qquad
\text{a.e.\ in $\Omega$.}
\]
Equivalently, there is $w^N_k \in \Lrm^\infty(\Omega;\R^m)$, namely $\norm{w^N_k}_{\Lrm^\infty} \leq C$ (uniformly in $k,N$), such that
\[
  - \Lcal^N_k u^N_k+\DD W_0(u^N_k) = f^N_k - w^N_k
\]
in a weak sense.
By Lemma~\ref{lem:apriori-space} below we thus get (recall that $p \in [2,\infty)$ from the definition of the external force $f$)
\begin{align}  \label{eq:apri_space}
  \norm{\nabla^2 u^N_k}_{\Lrm^p}\leq C (1+\norm{f^N_k}_{\Lrm^p}+\norm{f^N_k}_{\Lrm^2}^{q-1}).
\end{align}


The next lemma is an elliptic regularity result that is specifically taylored to our situation. We wish to point out that it allows for very general physically-motivated assumption on the elastic energy functional $\Wcal_0$, namely~\eqref{eq:Wgrowth}--\eqref{eq:Wmon}, but no structural assumptions like symmetry.

\begin{lemma} \label{lem:apriori-space}
With our usual assumptions from Section~\ref{sc:assume}, but excluding the mild convexity~\eqref{eq:muCP}, let $u\in (\Wrm^{1,2}_0 \cap \Lrm^q)(\Omega;\R^m)$ be any weak solution to
\begin{equation} \label{eq:space}
\left\{
\begin{aligned}
-\Lcal_t u +\DD W_0 (u)&=g \quad\text{ in $\Omega$,}\\
u|_{\partial\Omega}&=0.
\end{aligned} \right.
\end{equation}
If $g\in \Lrm^s(\Omega)$ for $s\in [2,\infty)$, then
\begin{equation} \label{eq:space_est}
 \norm{\nabla^2u}_{\Lrm^s}\leq C(1+ \norm{g}_{\Lrm^s}+\norm{g}_{\Lrm^{2}}^{q-1}).
\end{equation}
\end{lemma}

Note that the inhomogeneity with exponent $q-1$ is due to the $(q-1)$-growth of $\DD W_0$ via~\eqref{eq:DWgrowth}.

\begin{proof}
The existence of a solution $u\in (\Wrm^{1,2}_0 \cap \Lrm^q)(\Omega;\R^m)$ is guaranteed by the same variational argument as the one above. Analogous to~\eqref{eq:coerc} we find that
\begin{align}
\label{eq:coerc2}
\norm{\nabla u}_{\Lrm^2}\leq C(1+\norm{g}_{\Lrm^2}).
\end{align}
Alternatively, this can be deduced by using $u$ as a test function.

Next, we recall that the theory for elliptic operators implies that if
\[
-\Lcal_t u =g-\DD W_0 (u)\in \Lrm^s(\Omega)  \qquad\text{for some $s\in (1,\infty)$,}
\]
then
\begin{align} \label{eq:reg}
 \norm{\nabla^2u}_{\Lrm^s}\leq C \norm{\Lcal_t u}_{\Lrm^s}
\end{align}
since the coefficients are assumed to be Lipschitz continuous and the boundary of $\Omega$ is of regularity class $\Crm^{1,1}$.
For these results see~\cite{AgmonDouglisNirenberg59,AgmonDouglisNirenberg64} and also~\cite[Theorem 7.3]{GiaquintaMartinazzi12} (the standard scalar case is better-known and treated for instance in~\cite{Evans10book}).


Therefore, we are left to establish a bound on $\norm{\DD W_0(u)}_{\Lrm^s}$ for $s\in (1,\infty)$. In the case $d=2$ we find by Sobolev embedding that $\norm{u}_{\Lrm^s}\leq C\norm{u}_{\Wrm^{1,2}}$. Therefore,
 by~\eqref{eq:DWgrowth} and~\eqref{eq:coerc2},
\begin{align}
\label{eq:estw}
  \norm{\DD W_0(u)}_{\Lrm^s} \leq C(1+\norm{u}_{\Lrm^{s(q-1)}}^{q-1})
  \leq C(1+ \norm{g}_{\Lrm^{2}}^{q-1}).
\end{align}
Which implies, via~\eqref{eq:reg}, the wanted estimate.

In the following we will obtain the same estimate for $d=3$, which is assumed from now on until the end of the lemma. We will achieve this goal in several steps.

A first estimate concerns local $\Wrm^{2,2}$-regularity: For any ball $B_{5R} = B_{5R}(x_0) \subset\Omega$ ($x_0 \in \Omega$, $R > 0$) we will show
  \begin{align}
  \label{eq:local}
\int_{B_R} \abs{\nabla^2 u}^2 \dd x \leq \frac{C}{R^2} \biggl[ \int_{B_{5R}} \abs{\nabla u}^2 + \abs{g}^2\dd x \biggr].
\end{align}
  Since the system~\eqref{eq:space} is invariant under translation and scaling of coordinates, we may assume that $R=1$ and that the ball is centered around the origin, $x_0 = 0$. Indeed, $\tilde{u}(y)=u(R(y-z))$ solves~\eqref{eq:space} in the scaled and translated ball with $\tilde{g}(y)=R^2g(R(y-z)$ and $\tilde{W}_0:=R^2W_0$.

Thus, in the following we consider with no loss of generality that $u$ is a solution in the ball $B_5 = B_5(0)$. We take a cut-off function $\eta\in \Crm^1_0(B_2)$ such that $\eta \equiv 1$ on $B$ and $\supp \eta \subset 2B$, and further pick $h\in (0,1/2)$. We define the differential quotient in direction of the $k$'th unit vector $e_k$ by
\begin{align}
\label{eq:dq}
\DD^h_k g(x) := \frac{g(x+he_k)-g(x)}{h}= \dashint_0^h \partial_k g(x+se_k)\dd s.
\end{align}
Now, for $k\in \{1,..,d\}$ and $h \in (0,1/2)$ take $-\DD^{-h}_k(\eta^2 \DD^h_k(u)) \in \Wrm^{1,2}_0(B_3,\R^m)$ as a test function in~\eqref{eq:space} and employ partial summation to get
\begin{align*}
(I)+(II) &:= \int \sum_{i,j,\alpha,\beta} \DD^h_k (A^{\alpha,\beta}_{i,j}\partial_i u^\alpha) \partial_j(\eta^2 \DD^h_ku^\beta)\dd x+ \int \DD^h_k(\DD W_0(u))\cdot \DD^h_k u \, \eta^2 \dd x\\
&\phantom{:}= -\int g\cdot \DD^{-h}_k (\eta^2 \DD^h_k u)\dd x=:(III).
 \end{align*} 
 We begin with an estimate on $(I)$. By the product rule for difference quotients and~\eqref{eq:dq} we find
\begin{align*}
(I)&\geq   \int \sum_{i,j,\alpha,\beta} \DD^h_k(A^{\alpha,\beta}_{i,j}) \partial_i u^\alpha(\frarg + h\ee_k) \partial_j(\eta^2 \DD^h_ku^\beta)\dd x \\
&\qquad +2\int \sum_{i,j,\alpha,\beta} A^{\alpha,\beta}_{i,j}\DD^h_k(\partial_i u^\alpha) \eta (\partial_j \eta) \DD^h_ku^\beta\dd x\\
&\qquad +\int \sum_{i,j,\alpha,\beta} A^{\alpha,\beta}_{i,j}\DD^h_k(\partial_i u^\alpha) \eta^2 \DD^h_k (\partial_j u^\beta) \dd x\\
&\geq - C\norm{\DD A}_\infty \int \abs{\nabla u(\frarg + h\ee_k)}\abs{\DD^h_k\nabla u}\eta^2+ \abs{\nabla u} \abs{\DD^h_k u} \norm{\nabla \eta}_\infty \eta \dd x \\
&\qquad -{\norm{\nabla \eta}_\infty\norm{A}_{\infty}}\int \abs{\DD^h_k\nabla u}\abs{\DD^h_ku}\eta \dd x + \kappa \int \abs{\DD^h_k\nabla u}^2\eta^2 \dd x.
\end{align*}
Young's inequality then implies
\begin{align}
\label{eq:I}
(I)&\geq  
 \frac{\kappa}{2} \int \abs{\DD^h_k\nabla u}^2\eta^2 - C_{\kappa}(1+\norm{\nabla \eta}_\infty^2) \norm{A}_{\Crm^{0,1}}^2\int_{B_3} \abs{\nabla u}^2
+\abs{\DD^h_ku}^2 \dd x.
\end{align}
Next, we estimate $(II)$ by~\eqref{eq:Wmon} to find that  
\begin{align}
\label{eq:II}
(II)\geq -\mu\int\abs{\DD^h_k u}^2\eta^2 \dd x.
\end{align}
Finally, choosing $0<\eps<\kappa/2$, we find by Young's inequality and~\eqref{eq:dq} that
\begin{align} \label{eq:III}
(III) &\leq C_\eps \int (\abs{g}^2 +\norm{\nabla \eta}_\infty^2 \abs{\DD^h_k u}^2)\eta^2 \dd x \notag\\
&\qquad +\eps\dashint_0^{-h} \int \abs{\DD^h_k \partial_k u(x+se_k)}^2\eta^2\dd x \dd s. 
 \end{align} 
 Combining~\eqref{eq:I},~\eqref{eq:II} and~\eqref{eq:III}, we arrive at
 \begin{align*}
 \frac{\kappa}{2} \int_{B_1} \abs{\DD^h_k\nabla u}^2 \dd x
 &\leq C_{\kappa,\eps} 
 \int_{B_3} 
 \abs{\nabla u}^2 \dd x +C_\eps\int_{B_2}\abs{g}^2\dd x\\
 &\qquad + \eps\dashint_0^{-h} \int_{B_2} \abs{\DD^h_k \nabla u(x+se_k)}^2\dd x \dd s.
\end{align*}
To conclude the estimate, one wishes to absorb the $\eps$-term to the left hand side. This can be achieved via an interpolation result of Giaquinta--Modica type, see for instance~\cite[Lemma 13]{DieningEttwein08}. Thus, with a different constant $C_{\kappa,\eps}$,
 \begin{align*}
\int_{B_1} \abs{\DD^h_k\nabla u}^2 \dd x \leq C_{\kappa,\eps} \biggl[ \int_{B_{5}} \abs{\nabla u}^2 + \abs{g}^2\dd x \biggr].
\end{align*}
By letting $h\to 0$ we get the desired local estimate~\eqref{eq:local}. Moreover, embedding theory implies that
\begin{align}
\label{eq:local2}
\norm{u}_{\Lrm^s(B_R)}\leq C\norm{u}_{\Wrm^{2,2}(B_{R})}\leq \frac{C}{R} \bigl( \norm{u}_{\Wrm^{1,2}(B_{5R})}+\norm{g}_{\Lrm^2(B_{5R})} \bigr) 
\end{align}
for all $B_{5R}\subset\Omega$ and $s\in (1,\infty)$ if $d\in \{1,..,4\}$.

The next step is to get estimates near the boundary of $\Omega$. More precisely, we will show tangential differentiability up to the boundary, for which we will use a flattening argument. By the same scaling argument as before we assume the center point $0$ to be a boundary point, and that we have a one-to one diffeomorphism $\Psi \colon B_5 \to \R^d$, such that
\[
  \Psi(B_5\cap [\R^{d-1} \times \{0\}])=\partial\Omega\cap B_5
  \qquad\text{and}\qquad
  \Psi(B^+_5) = \Omega\cap B_5,
\]
where $B^+_R := B_R \cap [\R^{d-1} \times (0,\infty)]$.
By a straightforward transformation, we find that on the half-ball $B^+_5$ it holds that $\tilde{u}=u\circ\Psi:B^+_5\to\R^m$ is a weak solution to
\[
\widetilde{\Lcal}_t(\tilde u)+\DD W_0(\tilde u)=f\circ\Psi,
\] 
where 
 \[
 [\widetilde{\Lcal}_t]^\beta= \sum_j \partial_j\sum_{\alpha,i,k,l} J_{jk} A^{\alpha,\beta}_{kl}J_{li}\partial_iu^\alpha,  \qquad J_{ij}:=\partial_i\Psi^j.
 \]
We can assume that $\det \nabla\Psi>\kappa_1$ for some $\kappa_1>0$ depending on the prescribed boundary alone. Thus,~\eqref{eq:Lell} holds for $\kappa\kappa_1^2$. By~\eqref{eq:local} we get for the half balls that
   \begin{align}
  \label{eq:local3}
\dashint_{B_1^+} \abs{\partial_k \nabla \tilde{u}}^2 \dd x  \leq \frac{C}{R^2} \biggl[\dashint_{B_{5}^+} 
 \abs{\nabla \tilde{u}}^2\dd x+\dashint_{B_{5}^+}\abs{g}^2\dd x \biggr]
\end{align}
for $k\in\{1,...,d-1\}$. Up to this point the local estimates on the second derivatives are valid for any dimension $d\in \N$. In the following we will use our assumption $d=3$ and show that $u\in \Lrm^s(\Omega;\R^m)$ for all $s\in (1,\infty)$ and that
\begin{align}
\label{eq:ls}
\norm{u}_{\Lrm^s}\leq C(1+\norm{g}_{\Lrm^2}),
\end{align}
where $C$ depends on $s$, $\Omega$ and the constants of our assumptions.

%

%
%

For $z\in(0,1/2)$ we define
the function
$U^z(x,y):=\int_0^z \partial_z \tilde u(x,y,s)\dd s$.
We first observe that $U^z \in \Wrm^{1,2}(B_{1/2}\cap [\R^2 \times \{0\}])$. Indeed, 
\[
\int_{B_{1/2}\cap [\R^2 \times \{0\}]}\absBB{\partial_{x}\int_0^z \partial_z \tilde u(x,y,s)\dd s}^2\dd x\dd y
\leq \int_{B_{1/2}\cap [\R^2 \times \{0\}]}\bigg(\int_0^{1/2}\abs{\partial_x\partial_z\tilde u}\dd z\bigg)^2\dd x\dd y.
\]
 The last estimate holds for $\partial_y$ as well and is in both cases controlled by~\eqref{eq:local3}. 
 Sobolev embedding further implies
\[ 
  \norm{U^z}_{\Lrm^s(B_{1/2}\cap [\R^2 \times \{0\}])}\leq C\norm{U^z}_{\Wrm^{1,2}(B_{1/2}\cap [\R^2 \times \{0\}])} \qquad
  \text{for any $s\in (1,\infty)$.}
\]
From $\tilde u(x,y,z)=\int_0^z \partial_z \tilde u(x,y,t)\dd t$ we may therefore conclude
 \[
 \norm{\tilde{u}}_{\Lrm^s(B_{1/2}^+)}\leq C\norm{\tilde{u}}_{\Wrm^{1,2}(B_{5}^+)} + \norm{g}_{\Lrm^{2}(B_{5}^+)}.
 \]
To finish, we cover $\Omega$ with finitely many balls. For $x\in \overline{\Omega}$, there exists either $B_R(x)$, such that $B_{5R}(x)\subset \Omega$ or such that $B_{5R}(x)\cap\Omega$ is diffeomorphic to $B^+_{5R}(0)$. Since $\overline{\Omega}$ is compact we can choose a finite subfamily of balls for which either~\eqref{eq:local2} or~\eqref{eq:coerc2} holds.
This enables us to finish the proof for $d=3$ as in case $d=2$ and we get~\eqref{eq:space_est} via~\eqref{eq:reg} and~\eqref{eq:estw}.
\end{proof}

\subsection{Estimates for the time derivatives}

We test the $k$'th inequality of~\eqref{eq:ELweak} with $\xi=0$ and the $(k-1)$'th inequality with $u^N_k - u^N_{k-2}$, and divide by $h$ to find
\begin{align*}
 &\int_\Omega R(\delta_k^N)+R(\delta_{k-1}^N)+\nabla u^N_k : \Lbb^N_k :  \nabla\delta^N_k-\nabla u^N_{k-1} : \Lbb^N_{k-1} :  \nabla\delta^N_k
 \\
 &\qquad+(\DD W_0(u^N_k)-\DD W_0(u^N_{k-1}))\cdot\delta^N_k  -(f^N_k-f^N_{k-1})\delta^N_k \dd x - R(\delta^N_k+\delta^N_{k-1}) \dd x \leq 0.
\end{align*}
This can be transformed into
\begin{align*}
 &\int_\Omega R(\delta_k^N)+R(\delta_{k-1}^N)+(\nabla u^N_k-\nabla u^N_{k-1}) : \Lbb^N_k :  \nabla\delta^N_k+(\DD W_0(u^N_k)-\DD W_0(u^N_{k-1}))\cdot\delta^N_k   \\
&\qquad\,-(f^N_k-f^N_{k-1})\delta^N_k \dd x - R(\delta^N_k+\delta^N_{k-1}) - \nabla u^N_{k-1} : (\Lbb^N_k-\Lbb^N_{k-1}) :  \nabla\delta^N_k \dd x \leq 0.
\end{align*}
Divide by $h > 0$ and use the subadditivity of $R$ to get
\begin{align*}
  \int_\Omega \abs{\nabla\delta^N_k}^2
  &\leq \int_\Omega \frac{\abs{f^N_k-f^N_{k-1}}}{h} \cdot \abs{\delta^N_k}\dd x+\int_\Omega \frac{\abs{\Lbb^N_k-\Lbb^N_{k-1}}}{h}\cdot \abs{\nabla u^N_{k-1}}\cdot \abs{\delta^N_k}\dd x
  \\
  &\quad- \int_\Omega \frac{1}{h^2} \bigl[ \DD W_0(u^N_k) - \DD W_0(u^N_{k-1}) \bigr]  \cdot (u^N_k - u^N_{k-1}) \dd x .
\end{align*}
The first term on the right hand side we can estimate as
\begin{align*}
  \int_\Omega \frac{\abs{f^N_k-f^N_{k-1}}}{h} \cdot \abs{\delta^N_k} \dd x
  &\leq \biggl( \dashint_{t_{k-1}}^{t_k}\norm{\partial_t f}_{\Lrm^2} \dd t \biggr) \norm{\delta^N_k}_{\Lrm^2} \\
  &\leq C \biggl( \dashint_{t_{k-1}}^{t_k}\norm{\partial_t f}_{\Lrm^2} \dd t \biggr) \norm{\nabla \delta^N_k}_{\Lrm^2}.
\end{align*}
For the second term we use~\eqref{eq:Lcoeff}, the Cauchy--Schwarz and Poincar\'{e} inequalities as well as~\eqref{eq:space_est} to find
\begin{align*}
\int_\Omega \frac{\abs{\Lbb^N_k-\Lbb^N_{k-1}}}{h}\cdot \abs{\nabla u^N_{k-1}} \cdot \abs{\delta^N_k}\dd x
&\leq C \norm{\DD A}_\infty \norm{\nabla u^N_{k-1}}_{\Lrm^2}\norm{\nabla \delta^N_{k}}_{\Lrm^2} \\
&\leq C (1 + \norm{f^N_{k-1}}_{\Lrm^2}+\norm{f^N_k}_{\Lrm^2}^{q-1})\norm{\nabla \delta^N_{k}}_{\Lrm^2}.
\end{align*}
For the third term we use~\eqref{eq:Wmon} and the Poincar\'{e} inequality to get
\begin{align*}
  - \frac{1}{h^2} \int_\Omega \bigl[ \DD W_0(u^N_k) - \DD W_0(u^N_{k-1}) \bigr]  \cdot (u^N_k - u^N_{k-1}) \dd x
  &\leq \mu \int_\Omega \abs{\delta^N_k}^2 \dd x \\
  &\leq \mu C_P(\Omega)^2 \int_\Omega \abs{\nabla \delta^N_k}^2 \dd x,
\end{align*}
where we recall that by $C_P(\Omega) > 0$ we denote the Poincar\'{e} constant of $\Omega$.
Hence, combining, we get
\begin{equation} \label{eq:apri_time}
 \norm{\nabla \delta^N_k}_{\Lrm^2} \leq \frac{C}{1-\mu C_P(\Omega)^2} \biggl[ 1 + \dashint_{t_{k-1}}^{t_k}\norm{\partial_t f}_{\Lrm^2} \dd t +\norm{f^N_k}_{\Lrm^2}+\norm{f^N_k}_{\Lrm^2}^{q-1} \biggr] 
\end{equation}
By~\eqref{eq:muCP}, the constant on the right is greater than zero.


%

%

\subsection{H\"older continuity of the gradient}
\label{ssec:Hgrad}
In this section only we additionally assume that $p > d$ (see the statement of Theorem~\ref{thm:main}).

For Borel subsets $E\subset\R^d$ with positive and finite Lebesgue measure we will use the notation 
\[
\mean{f}_E:=\dashint_E f \dd x =\frac{1}{\abs{E}}\int_Ef \dd x.
\]
We also define
\[
  u^N(t):=\frac{t-t_{k-1}}{h}u^N_k+\frac{t_k-t}{h}u^N_{k-1}  \qquad\text{for $t \in (t_{k-1},t_k]$,\; $k = 1,\ldots,N$.}
\]

For $p\in [2,\infty)$, $a\in (1,\infty]$, and $f\in \Wrm^{1,a}(0,T;\Lrm^p(\Omega,\R^m))$ we find by embedding that $f\in \Crm^{0,(a-1)/a}([0,T];\Lrm^p(\Omega,\R^m))$. Therefore,~\eqref{eq:apri_space} and~\eqref{eq:apri_time} imply
\begin{align}
\label{eq:apri_uN}
  \norm{\nabla^2 u^N}_{\Lrm^\infty(\Lrm^p)}
  + \norm{\nabla \dot{u}^N}_{\Lrm^a(\Lrm^2)} \leq C,
\end{align}
uniformly in $N$. Thus,
\begin{align}
\label{eq:infty2p}
  u^N \in \Lrm^\infty(0,T;\Wrm^{2,p}(\Omega;\R^m)),  \qquad
  \dot{u}^N \in \Lrm^a(0,T;\Wrm^{1,2}(\Omega;\R^m)).
\end{align}
and the respective norms are uniformly bounded.

Recall that $\Wrm^{2,p}(\Omega;\R^m) \embed \Crm^{1,\alpha}(\cl{\Omega};\R^m)$ for some $\alpha \in (0,1)$ if $p\in (d,\infty)$. Thus, $u^N \in \Lrm^\infty(0,T;\Crm^{1,\alpha}(\cl{\Omega};\R^m))$ and in fact the embedding is compact into $\Lrm^r(0,T;\Crm^{1,\alpha}(\cl{\Omega};\R^m))$ for any $r \in [1,\infty)$ and a any smaller $\alpha \in (0,1)$.
We will show that if $a\in (1,\infty]$ and $p\in (d,\infty)$, then $\nabla u^N$ is uniformly Lipschitz continuous with respect to the metric
\[
  \rho \bigl( (t,x),(s,y) \bigr) := \abs{t-s}^\zeta +\abs{x-y}^\alpha,
\]
for any $\alpha \in (0,1)$ such that $u^N \in \Lrm^\infty(\Crm^{1,\alpha})$, and where
\[
  \zeta = \frac{\alpha(a-1)}{\big(\frac{d}{2}+\alpha\big)a} = \frac{\alpha}{b}, \qquad
  b=\Big(\frac{d}{2}+\alpha\Big)\frac{a}{a-1}.
\]

First note that because of the zero Dirichlet boundary values, we can always extend $u^N$ onto the whole space. By Campanato's integral characterization of H\"{o}lder continuity~\cite{Campanato63} (also see Section~III.1 in~\cite{Giaquinta83book} or Section~2.3 in~\cite{Giusti03book}), we need to show
\[
  \dashint_{t}^{t+r^b}\!\!\dashint_{B_r(x)}\abs{\nabla u^N-\mean{\nabla u^N}_{(t,t+r^b)\times B_r(x)}} \dd x\dd t
  \leq Cr^\alpha
\]
for all $(t,t+r^{b})\times B_r(x)\subset [0,T]\times \R^d$, $r > 0$. Indeed, one can check easily that this \enquote{parabolic} version follows from the usual one via the transformation $g(s,x) = bs^{b-1}f(s,x)$.

 For easier reading we assume $(t,x)=(0,0)$ and estimate
\begin{align*}
 &\dashint_0^{r^{b}}\!\!\dashint_{B_r}\abs{\nabla u^N-\mean{\nabla u^N}_{(0,r^{b})\times B_r}} \dd x\dd t\\
&\qquad\leq \sup_{0 \leq t \leq r^{b}} \; \dashint_{B_r}\abs{\nabla u^N(t,x)-\mean{\nabla u^N(t)}_{B_r}} \dd x \\
&\qquad\qquad + \dashint_0^{r^{b}}\abs{\mean{\nabla u^N}_{(0,r^{b})\times B_r}-\mean{\nabla u^N(t)}_{B_r}} \dd t\\
&\qquad\leq Cr^\alpha +C \underbrace{\dashint_0^{r^{b}}r^{b}\dashint_{B_r}\abs{\nabla \dot{u}^N} \dd x \dd t}_{=: (I)} 
\end{align*}
Here we used the a-priori $\Lrm^\infty(\Crm^{1,\alpha})$-regularity on the first integral and the Poincar\'{e} inequality in the time direction on the second integral.

To bound~$(I)$, use the $\Lrm^a(\Lrm^2)$ estimate of $\nabla \dot{u}^N$ and H\"older's inequality to get with $a' = a/(a-1)$ that
\begin{align*}
 (I)&\leq\int_0^{r^{b}} \Bigg(\dashint_{B_r}\abs{\nabla \dot{u}^N}^2 \dd x\bigg)^{1/2}\dd t\\
&\leq r^{b/a'}\bigg(\int_0^{r^{b}}\Bigg(\dashint_{B_r}\abs{\nabla \dot{u}^N}^2 \dd x\bigg)^{a/2}\dd t\bigg)^{1/a}\\
&= Cr^{b/a'-d/2}\norm{\nabla \dot{u}^N}_{\Lrm^a(\Lrm^2)} \\
&\leq Cr^\alpha.
\end{align*}

Finally, the $d$-H\"{o}lder continuity also implies $\zeta$-H\"{o}lder continuity jointly in space and time. To see this, we can estimate, since $\zeta< \alpha$,
\begin{align*}
  \abs{t-s}^{\zeta} + \abs{x-y}^\alpha &\leq 2 \max \bigl\{ (\abs{t-s} + \abs{x-y})^{\zeta}, (\abs{t-s} + \abs{x-y})^\alpha \bigr\} \\
  &\leq 2(1+(T + \diam(\Omega))^{\alpha-\zeta}) (\abs{t-s} + \abs{x-y})^{\zeta},
\end{align*}
where we have to consider the cases $\abs{t-s} + \abs{x-y} \leq 1$ and $\abs{t-s} + \abs{x-y} > 1$, separately.

\subsection{H\"older continuity of the solution}
\label{ssec:Hsol}

By a similar argument to the one in the last section, we will show that $u^N$ is uniformly H\"older continuous. We only need to to consider the case $p\in [2,d]$, since otherwise the uniform H\"older continuity of $u^N$ follows from~\eqref{eq:infty2p}.


Recall first (see~\eqref{eq:infty2p}) that $\nabla^2 u^N \in \Lrm^\infty(\Lrm^p)$ implies that $\nabla u^N \in \Lrm^\infty(\Wrm^{1,s}_0)$  for all $s\in [1,\frac{pd}{d-p}]$ in case $p<d$ and all $s\in (1,\infty)$ in case $d=p$. Since $\frac{pd}{d-p}\geq 6$, there exists an $\alpha\in (0,1)$ such that $u^N\in \Lrm^\infty(\Crm^{0,\alpha})$. Furthermore, as $\nabla \dot{u}^N\in \Lrm^a(\Lrm^2)$, we find that $\dot{u}^N\in \Lrm^a(\Lrm^s)$, for $s\in [1,6]$ if $d=3$ and $s\in (1,\infty)$ for $d=2$. By~\eqref{eq:infty2p},
\[
  u^N \in \Lrm^\infty(0,T;\Crm^{0,\alpha}_0(\Omega,\R^m)) \cup \Wrm^{1,a}(0,T;\Lrm^{s}(\Omega,\R^m))
\]
and the norms can be correspondingly estimated by an $N$-independent constant.

We can therefore argue exactly as before. Indeed, we will show that $u^N$ is Lipschitz continuous with respect to the metric
\[
  \tilde{\rho} \bigl( (t,x),(s,y) \bigr) := \abs{t-s}^\gamma +\abs{x-y}^\alpha,
\]
for any $\alpha \in (0,1)$ such that $u^N \in \Lrm^\infty(\Crm^{0,\alpha})$ (our $\alpha$ here is different from the one in the previous section), and where
\[
  \gamma = \frac{\alpha(a-1)}{\big(\frac{d}{2}+\alpha\big)a} = \frac{\alpha}{b}, \qquad
  b=\Big(\frac{d}{2}+\alpha\Big)\frac{a}{a-1}.
\]
Again we use the variables $a'=\frac{a}{a-1}$. By Campanato's integral characterization of H\"{o}lder continuity, we need to show
\[
  \dashint_{t}^{t+r^b}\!\!\dashint_{B_r(x)}\abs{\nabla u^N-\mean{\nabla u^N}_{(t,t+r^b)\times B_r(x)}} \dd x\dd t
  \leq Cr^\alpha
\]
for all $(t,t+r^{b})\times B_r(x)\subset [0,T]\times \R^d$, $r > 0$. 

We assume $(t,x)=(0,0)$ and estimate as before and additionally with Poincar\'{e}{}'s inequality,
  \begin{align*}
 &\dashint_0^{r^{b}}\!\!\dashint_{B_r}\abs{u^N-\mean{u^N}_{(0,r^{b})\times B_r}} \dd x\dd t\\
&\qquad \leq Cr^\alpha+ C \dashint_0^{r^{b}}r^{b}\dashint_{B_r}\abs{\dot{u}^N} \dd x \dd t \\
&\qquad \leq Cr^\alpha+ C r^\frac{b}{a'}\int_0^{r^{b}}\bigg(\dashint_{B_r}\abs{\nabla \dot{u}^N}^2 \dd x\bigg)^\frac{a}{2} \dd t \\
&\qquad \leq Cr^\alpha. 
\end{align*}
Hence, we find that $u^N$ is uniformly $\gamma$-H\"older continuous.
\begin{remark}
\label{rem:H}
The estimates in the last two sections are of general nature. Indeed, what is shown here is
\[
[u]_{\Crm^{0,\gamma}([0,T]\times\overline{\Omega})}\leq C\Big(\norm{u}_{\Lrm^\infty(\Wrm^{2,p})}+\norm{\dot{u}}_{\Lrm^a(\Wrm^{1,2})}\Big),
\]
where we can choose $\alpha\in (0,\min\{1,(2p-d)/p\})$ and $\gamma$ accordingly.

Moreover, if $p > d$, then
\[
[\nabla u]_{\Crm^{0,\zeta}([0,T]\times\overline{\Omega})}\leq C\Big(\norm{\nabla u}_{\Lrm^\infty(\Wrm^{1,p})}+\norm{\nabla\dot{u}}_{\Lrm^a(\Lrm^{2})}\Big)
\]
for $\alpha\in (0,\min\{1,(p-d)/p\})$ and $\zeta$ chosen accordingly.
\end{remark}

\section{Passing to the limit and proof of Theorem~\ref{thm:main}}  \label{sc:proof}

The a-priori estimate~\eqref{eq:apri_uN} implies that there exists a (non-relabelled) subsequence such that
\[
  u^N \toweakstar  u  \quad\text{in $\Lrm^\infty(\Wrm^{2,p})$} \qquad\text{and}\qquad
  u^N \toweak  u  \quad\text{in $\Wrm^{1,a}(\Wrm^{1,2})$}.
\]
By the weak compactness in reflexive Banach spaces we can furthermore assume
\[
  u^N\toweak u  \quad\text{in $\Lrm^r(\Wrm^{2,p})$}
   \qquad  \text{for any $r\in(1,\infty)$.}
\]

We rewrite~\eqref{eq:discrete1} into a continuous form. Observe that on $(t_{k-1},t_k]$ we have $\dot{u}^N=\delta^N_k$ and
\[
 \nabla u^N(t) = \frac{t-t_{k-1}}{h}\nabla u^N_k+\frac{t_k-t}{h}\nabla u^N_{k-1}= \nabla u^N_k +\frac{t_k-t}{h}\nabla (u^N_{k-1}-u^N_{k}).
\]
We also set 
\[
  k_N(t) := \setb{k\in\{1,\ldots,N\}}{t\in (t_{k-1},t_k]}  \qquad
  \text{for $t\in [0,T]$.}
\] 
Therefore,~\eqref{eq:discrete1} here reads as
\begin{align}
&\int_0^T\!\!\int_\Omega R(\dot{u}^N(t)) - \nabla u^N_{k_N(t)} : \Lbb^N_{k_N(t)} : \nabla(\xi(t)-\dot{u}^N(t)) \dd x \dd t\notag \\
&\qquad\qquad + \int_0^T\!\!\int_\Omega \bigl[-\DD W_0(u^N_{k_N(t)})+f^N_{k_N(t)} \bigr]\cdot(\xi(t)-\dot{u}^N(t)) \dd x \dd t\notag \\
&\qquad \leq \int_0^T\!\!\int_\Omega R(\xi) \dd x\dd t  \label{eq:tolimit}
\end{align}
for all $\xi \in \Lrm^1(0,T;\Wrm^{1,2}_0(\Omega;\R^m))$ (first use only $\xi$ that are piecewise constant with respect to $\{t^N_0,t^N_1,\ldots,t^N_N\}$ and then argue by density).
 
Using the H\"older continuity of $u$, we find by the Arzel\`{a}--Ascoli theorem a subsequence such that
\[
  u^N \to u
  \quad\text{in $\Crm^{0,\gamma}([0,T]\times\overline{\Omega})$}\qquad
  \text{for $0<\gamma< \frac{\alpha(a-1)}{\big(\frac{d}{2}+\alpha\big)a}$}\qquad
  \text{if $p\in [2,d]$.}
\]
and 
\[
  u^N \to u
  \quad\text{in $\Crm^{1,\zeta}([0,T]\times\overline{\Omega})$}\qquad
  \text{for $0<\zeta< \frac{\alpha(a-1)}{\big(\frac{d}{2}+\alpha\big)a}$}\qquad
  \text{if $p\in (d,\infty]$.}
\]
Here, $\alpha$ is defined as above via the respective Sobolev embedding in space, see Remark~\ref{rem:H}. By the equi-continuity we also know that in both cases
\[
  u^N_{k_N(t)}(x)\to  u(t,x)  \quad\text{in $\Crm^{0,\gamma}(\overline{\Omega})$} \qquad
  \text{for every $t \in (0,T]$ and any $0<\gamma<\frac{\alpha(a-1)}{\big(\frac{d}{2}+\alpha\big)a}$.}
\]


By the convexity and lower semicontinuity of $R$ as well as the assumptions on $\DD W_0$ (continuity) and $f \in \Crm^0(\Lrm^2)$, we get
\begin{align*}
&\int_0^T\!\!\int_\Omega R(\dot{u})+\bigl[-\DD W_0(u)+f\bigr]\cdot(\xi-\dot{u}) \dd x \dd t\\
&\qquad \leq \liminf_{N\to\infty}\int_0^T\!\!\int_\Omega R(\dot{u}^N)+ \bigl[-\DD W_0(u^N_{k_N})+f^N_{k_N} \bigr]\cdot(\xi-\dot{u}^N) \dd x \dd t 
\end{align*}
for all $\xi\in \Lrm^1(0,T;\Wrm^{1,2}_0(\Omega;\R^m))$.

The term of the regularizer needs special attention. Rellich's compactness theorem implies that $\Lrm^\infty(\Wrm^{2,p}) \cap \Wrm^{1,a}(\Wrm^{1,2}_0)$ is compactly embedded in $\Crm^{\beta}(\Wrm^{1,2})$ for some $\beta>0$, see~\cite{Simon87}. Therefore by passing to yet another subsequence, we find that $\nabla u^N\to\nabla u$ in the strong topology of $\Crm^{\beta}(\Wrm^{1,2})$, in particular
\[
  \norm{\nabla u^N(t) - \nabla u(t)}_{\Lrm^2} \to 0 \quad\text{uniformly in $t \in [0,T]$.}
\]

Consequently,
\[
\int_0^T\!\!\int_\Omega\nabla u^N_{k_N} : \Lbb^N_{k_N} : (\xi-\nabla(\dot{u}^N)\dd x\dd t\to \int_0^T\!\!\int_\Omega\nabla u: \Lbb : \nabla(\xi-\dot{u})\dd x\dd t.
\]
Hence, letting $N\to\infty$ in~\eqref{eq:tolimit}, we get
\begin{align*}
&\int_0^T\!\!\int_\Omega R(\dot{u}) - \nabla u : \Lbb : \nabla(\xi-\dot{u})+\bigl[-\DD W_0(u)+f\bigr]\cdot(\xi-\dot{u}) \dd x\dd t \\
&\qquad \leq \int_0^T\!\!\int_\Omega R(\xi) \dd x\dd t,
\end{align*}
for all $\xi \in \Lrm^1(0,T;\Wrm^{1,2}_0(\Omega;\R^m))$. Therefore, the limit inequality~\eqref{eq:strongsol_ineq} is established and our $u$ is indeed a strong solution to~\eqref{eq:PDE}.
Hence, combining all of the above assertions, Theorem~\ref{thm:main} is proved.




\providecommand{\bysame}{\leavevmode\hbox to3em{\hrulefill}\thinspace}
\providecommand{\MR}{\relax\ifhmode\unskip\space\fi MR }
\providecommand{\MRhref}[2]{%
  \href{http://www.ams.org/mathscinet-getitem?mr=#1}{#2}
}
\providecommand{\href}[2]{#2}

\end{document}